\documentclass{article} 
\usepackage{amsfonts}

\newtheorem{Th}{Theorem}

\newtheorem{Cor}{Corollary}
\newtheorem{Def}{Definition}
\newtheorem{Rem}{Remark}

\newcommand{\eqdef}{\stackrel{{\rm def}}{=}}

\newcommand{\Id}{\mbox{\rm Id}}

\begin{document}
\title{ Closed manifolds admitting
metrics with the same geodesics}
\author{Vladimir S. Matveev\thanks{Mathematisches Institut, Universit\"at 
Freiburg, 79104 Germany \ 
 matveev@email.mathematik.uni-freiburg.de}}
\date{}
\maketitle

\begin{abstract} { 
The goal of this survey is to give a list of resent results about  topology of manifolds admitting different metrics with the same geodesics. We emphasize the 
role of the theory of integrable systems in obtaining  these results. 
}
\end{abstract}

\date{}

\section{Introduction} 
\subsection{Definitions}
\begin{Def}
Let $g$ be a Riemannian  metric on a manifold  $M^n$ of dimension $n\ge 2$. A Riemannian 
metric $\bar g$ on  $M^n$ is called 
{\bf  geodesically equivalent} to $g$, if every geodesic of $\bar g$, 
considered 
as an unparameterized curve, is a geodesic of $g$. 
\end{Def}

Trivial examples of geodesically equivalent metrics can be obtained by 
considering proportional metrics 
 $g$ and $C\cdot  g$, where $C$ is a positive 
constant.

\begin{Def}  A manifold $M^n$   
is called {\bf geodesically rigid},    if every  two   geodesically equivalent 
Riemannian  metrics  on $M^n$   are proportional.   
\end{Def}

In other words, on 
geodesically rigid manifolds, unparameterized geodesics define the metric (modulo multiplication by a constant).

\subsection{ History}

The theory of  geodesically equivalent metrics has a  long and  fascinating 
history that goes back to the works of  Beltrami, Dini and Levi-Civita. 

Beltrami~\cite{Beltrami} 
was the first to  observe that two nonproportional 
 metrics (even on closed  manifolds) can have the 
same geodesics. 

At  the end of his paper~\cite{Beltrami}, Beltrami formulated the  problem of 
describing all geodesically equivalent metrics (for  surfaces.)
 It is not clear  from the text whether he assumed  the  
  local or the   global description; actually, his motivation came from 
a certain  problem of cartography,  which requires the  global  setting.
Nevertheless, partially  because of strong results of     
   Dini,  Levi-Civita, Weyl, 
E. Cartan
  and  Eisenhart, 
  the theory of geodesically equivalent 
metrics was mostly a local geometry.

Locally, 
  in a neighborhood of almost every  point,  
 a complete  
description of geodesically equivalent metrics 
has been    given by  Dini~\cite{Dini} 
for surfaces and Levi-Civita~\cite{Levi-Civita} 
for manifolds of arbitrary dimension.  As a corollary of this description, one can show that, at least for  dimensions two and three, 
   every open  manifold has  non-proportional  geodesically  
 equivalent metrics.

Later, geodesically equivalent metrics were considered by 
Weyl, E. Cartan and  Eisenhart. Weyl studied
geodesically equivalent metrics on the tensor level and found a 
few tensor reformulations of geodesic equivalence.  
 One of his  most
remarkable results  is  
  the  construction~\cite{Weyl2} of the projective Weyl tensor: 
 if two metrics are geodesically equivalent, then their  
  projective Weyl tensors coincide. 
E. Cartan~\cite{cartan}  studied  geodesic equivalence on the 
level of   affine connections. He introduced the so-called projective 
connection, which allows   reconstruction of   
 geodesics as unparameterized curves.  
In his book~\cite{Eisenhart},  Eisenhart   systematically  
applied both methods and   obtained   a series of local  results.

It is clear that the classics such as Lie~\cite{Lie},
 Painlev\'e~\cite{painleve}, 
  Levi-Civita and Eisenhart 
understood  well the connection between integrable 
systems and  geodesically equivalent metrics.  But they did not use it, 
probably because  they mostly were interested in the local aspects  of
geodesically equivalent metrics. 

 Global aspects  have  been intensively  
studied since 50th, firstly by 
French (Lichnerowicz), Soviet (Rashevskii and Solodovnikov)   and 
Japanese  (Yano)   geometry  schools.  
But, probably because of 
 the influence of  earlier researcher,  all known global 
 results  require fairly  strong additional  geometrical assumptions.

 Roughly  speaking, one   takes  some  geometric assumption written in 
  tensor form, combines  it  with one 
of the tensor reformulations of geodesic equivalence  and   deduces  some new 
  object   with global  geometric properties, 
   see the   surveys of Mikes~\cite{Mikes2} and Aminova~\cite{aminova}.

Recently, it was found that integrable systems provides a very 
effective tool for understanding the topology of the manifolds admitting 
geodesically equivalent metrics. We describe the connection 
between geodesically equivalent metrics and integrable geodesic flows in Section~\ref{int}.   Roughly speaking, the existence of $\bar g$  geodesically equivalent to $g$ allows one to construct  integrals for the geodesic flow 
of $g$, see Theorem~\ref{integrability} for details.  If the metric $g$ and 
$\bar g$ are strictly-non-proportional, the geodesic flow of $g$ is 
Liouville-integrable.

\section{Resent results}
\begin{Th}[\cite{hyperbolic}] \label{main}
Let $M^n$ be a closed connected   manifold. Suppose  two
 non-proportional  Riemannian metrics $g$, $\bar g$ on $M^n$  
are   geodesically equivalent.
If  the  fundamental group   of $M^n$ is infinite, then 
  there exist  $r\in \{1,2,...,n-1\}$,   a Riemannian 
metric $\tilde g$  and
  foliations 
$B_r$ (of dimension $r$) and $B_{n-r}$ (of dimension $n-r$) such that,  
in a   neighborhood $U(p)$  of every  point $p\in M^n$, there exist  
coordinates 
$$
(\bar x,\bar y)= \bigr((x_1,x_2,...,x_r),(y_{r+1},y_{r+2},...,y_n)\bigl)
$$
 such that  the  $x$-coordinates are 
constant on every fiber of the foliation $B_{n-r}\cap U(p)$, 
the $y$-coordinates are 
constant on every fiber of the foliation $B_{r}\cap U(p)$,   and 
 the metric  $\tilde g$ has the  block-diagonal form   
\begin{equation} \label{product}
ds^2= \sum_{i,j=1}^r G_{ij}(\bar x)dx_idx_j + \sum_{i,j=r+1}^{n} 
H_{ij}(\bar y)dy_idy_j, 
\end{equation} 
where  the first  block depends on the first $r$ coordinates and 
the second  block depends on the remaining  $n-r$ coordinates. 
\end{Th}

Theorem~\ref{main}  already gives us the complete list of all geodesically
 rigid closed   surfaces~\cite{ERA}: a closed connected  surface is geodesically
 rigid if and only if its Euler characteristic is negative.  

 More precisely, a closed  connected surface of negative Euler characteristic 
 has infinite  fundamental group, and admits no one-dimensional foliation,
 so it is geodesically rigid. 

A closed connected surface of  nonnegative Euler characteristic is  not 
geodesically rigid: it is diffeomorphic  to   the
 sphere or  the projective plane or the torus or the Klein bottle. 
Examples of nonproportional  geodesically  
 equivalent metrics  on  the
 sphere and on   the projective plane were essentially constructed by 
Beltrami~\cite{Beltrami}. 
Since the geodesics of every   flat
 metric are straight lines, every  two 
flat metrics on the torus (or  on the Klein bottle) are geodesically related 
 (that is,   there  exists a diffeomorphism that 
takes the geodesics of the first metric  to the geodesics of the  second), so 
the torus or the Klein bottle are not geodesically rigid as well. 

For dimension three, a direct corollary~\cite{hyperbolic} of Theorem~\ref{main} is 

\begin{Cor}\label{main1}
Suppose $M^3$ is a connected closed manifold. 
Suppose  there exist  nonproportional Riemannian metrics  
on $M^3$  that are   geodesically 
equivalent. Then, modulo the  Poincar\'e conjecture, 
 $M^3$ is finitely covered  by the sphere $S^3$ or by the product 
$F^2\times S^1$, where $F^2$ is a closed surface.
\end{Cor}

It appears that  Corollary~\ref{main1} is true~\cite{short,threemanifold} also 
without assuming the  Poincar\'e conjecture:

\begin{Th} \label{dim3}
Let nonproportional   Riemannian metrics $g$ and $\bar g$ be 
  geodesically equivalent on 
 a  closed connected  
three-dimensional manifold  $M^3$.  Then the manifold 
 is homeomorphic either to a lens space or to a Seifert manifold with 
zero Euler number. Every  lens space and every Seifert manifold with 
zero Euler number admits   geodesically  equivalent metrics which are nonproportional. 
\end{Th}  

Theorems~\ref{main},\ref{dim3} give us 
  a complete list of closed connected manifolds of dimension two and three
  admitting nontrivial geodesic  equivalence. 
It  will  be much more complicated to obtain 
  such list in every dimension.
But   still, in every dimension $n\ge 2$,  there exists infinitely many 
 geodesically 
rigid manifolds~\cite{starrheit,hyperbolic}: 
\begin{Th} \label{hyperbolic} 
Every  closed connected   manifold  admitting a Riemannian metric of 
negative sectional curvature  is  geodesically rigid. 
\end{Th} 
Note that in view of result of Borel~\cite{borel}, in every dimension there 
exist infinitely many closed  manifolds admitting   metrics
 of negative sectional curvature.

\section{Methods and ideas} 
\subsection{Integrability for  the geodesic flows of 
 geodesically equivalent metrics}\label{int} 
New methods for the global (= on closed or complete manifolds)  
investigation of geodesically equivalent metrics are based on the 
following observation~\cite{MT,ERA,dedicata}:  
  the existence of  $\bar g$ geodesically equivalent to
 $g$ allows one to construct commuting  integrals for the geodesic flow of 
$g$. 

Let $g=(g_{ij})$ and  $\bar g=(\bar g_{ij})$ be Riemannian metrics on 
a manifold $M^n$. 
Consider the (1,1)-tensor  $L$ 
given by the formula   
\begin{eqnarray}
L^i_j &\eqdef & \left(\frac{\det(\bar g)}{\det(g)} \right)^\frac{1}{n+1} \bar g^{i\alpha} g_{\alpha j}. \label{l}
\end{eqnarray}
Then, $L$ determines the family $S_t$, $t\in R$, of $(1,1)$-tensors 
\begin{equation}\label{st}
 S_t\eqdef \det(L - t\ \mbox{\rm Id})\left(L-t\ \mbox{\rm Id}\right)^{-1}. 
 \end{equation}
\begin{Rem} 
Although $\left(L-t\ \Id\right)^{-1}$ is not defined for 
$t$ 
lying 
in the spectrum of $L$, the tensor  $S_t$  is well-defined 
for every  $t$.  Moreover,   $S_t$
is a  polynomial
 in $t$ of degree $n-1$ 
with coefficients being  (1,1)-tensors.  
\end{Rem}
\noindent We will identify the tangent and  cotangent bundles of $M^n$ by $g$. 
This identification allows us to transfer the 
 natural  Poisson structure from  $T^*M^n$ to   $TM^n$.

\begin{Th}\label{integrability}
 If $g$, $\bar g$ are geodesically 
 equivalent,   
then, for every  $t_1,t_2\in R$, the functions 
\begin{equation}\label{integral}
I_{t_i}:TM^n\to R, \ \ I_{t_i}(v)\eqdef g(S_{t_i}(v),v)
\end{equation}
are commuting integrals for the geodesic flow
 of  $g$. 
\end{Th}

In other direction these  theorem is wrong; a counterexample could be found in~\cite{quantum}. 

\begin{Th}[\cite{quantum}]\label{integrability1}
Suppose for    
  every  $t\in R$ the function 
$
I_{t}$ given by (\ref{integral})
is an  integral for the geodesic flow
 of  $g$. If the Nijenhuis torsion $N_L$ vanishes, the metrics are 
geodesically equivalent.
\end{Th}

\begin{Th}[\cite{benenti}]\label{nijenhuis}
Let $g$, $\bar g$ be  geodesically 
 equivalent. Then the Nijenhuis torsion $N_L$ vanishes. 
\end{Th}

\subsection{What is special in these integrals?}

\begin{Def} Two metrics $g$ and $\bar g$ 
 are {\bf strictly-non-proportional} at $x\in M^n$, if  all 
roots of $P(t):=det(g- t\bar g)$ are simple.
\end{Def}

Let us assume that there exists a point where the 
metrics are strictly-non-proportional. 
Then, it is so~\cite{hyperbolic}  at almost every point, 
and the family $I_t$ contains $n$ integrals that are functionally independent 
almost everywhere. Hence, the geodesic flow of the metric is 
Liouville-integrable.

Let us note that 
\begin{itemize} 
\item the integrals are quadratic in velocities, 
\item at every point $x\in M^n$, the integrals (considered as quadratic forms) 
 can be simultaneously diagonalizable.  
\end{itemize} 
Integrable systems  with such properties~\cite{stackel} 
are known as St\"ackel systems.  Locally, in a given coordinate system, 
it can be defined by using a ($n\times n$)-matrix such that 
its columns  depend on the corresponding coordinate only. 

Not every st\"ackel system can come from geodesically equivalent metrics. 
The additional assumption is that the st\"ackel matrix can be chosen to be  
 a  Vandermonde matrix. 

St\"ackel systems with such property were also intensively studied. 
One of the reasons for it that they satisfy Robertson's~\cite{robertson} condition, which 
imply that its quantization is quantum-integrable as well. The second reasons 
is that  all st\"ackel systems coming from physics are of this type (or, a 
degeneration of systems of this type).

It is possible to show~\cite{quantum}, that this extra-condition 
is also a sufficient condition for the  existence of geodesically 
equivalent metrics.

The systems with this condition appear independently and under different
 names (L-systems, Benenti systems, quasi-bihamiltonian systems) in works 
of  different authors~\cite{Benenti1,Benenti2,Benenti3,benent,I,crampin}.

\subsection{If the metrics are strictly non-proportional} 

 Theorem~\ref{integrability}  can be used most  efficiently
when there exists a point of a manifold  where the 
 metrics are strictly-non-proportional.
   Then,  the geodesic flow of $g$ is Liouville-integrable, and we can apply 
the well-developed machinery of integrable systems. For example, 
 the following theorem follows directly  from  Theorem~\ref{integrability}
 and Taimanov~\cite{Taimanov}. 

\begin{Th}\label{dubrovin} 
Suppose $M^n$ is a connected closed manifold. 
Let the  real-analytic Riemannian metrics  
$g$  and $\bar g$ on $M^n$ 
 be geodesically  
equivalent. Suppose  there exists a point of the  manifold  where the metrics are strictly-non-proportional. 
Then, the following statements hold. 
\begin{enumerate} 
\item 
 The first Betti number $b_1(M^n)$ is 
not  greater than $n$.

\item The fundamental group $\pi_1(M^n)$ is virtually Abelian. 

\end{enumerate}
\end{Th}

The integrals are quadratic in velocities. Combining this fact with 
 topological obstructions~\cite{Kol}  for the existence of   
 quadratically-integrable geodesic flows   on  closed  surfaces,
 we obtain Theorem~\ref{hyperbolic}  for dimension two~\cite{ERA}.

Note that, in view of results~\cite{Waldhausen},
 Corollary~\ref{main1} follows from 
Theorem~\ref{dubrovin} under the additional assumption that  the metrics are  
real-analytic and that there exists a point where the 
metrics are strictly-non-proportional .

\subsection{Geodesic equivalence and zero entropy} \label{kruglikov}
All results of this section are joint with Kruglikov; 
the proofs will be published elsewhere.

It is expected, that an integrable geodesic flow has zero topological 
 entropy. 
This is not always the case.  There are examples~\cite{BT1,BT2} of integrable flows 
with non-zero entropy. But still it is possible to show  that
 if the integrals come from strictly-non-proportional geodesically equivalent metrics by applying  
Theorem~\ref{integrability}, the topological entropy of the geodesic flow 
must be zero.

\begin{Th}  
Suppose the Riemannian metrics $g$, $\bar g$ on
a closed connected $M^n$ are geodesically equivalent.
 Suppose there exists a point where the metrics are strictly non-proportional.
 Then, the topological entropy of the geodesic flow of $g$ is zero.  
\end{Th} 
Combining this theorem with the famous Yomdin's Theorem~\cite{gromov}, 
we obtain

\begin{Cor} 
Suppose the Riemannian metrics $g$, $\bar g$ on
a closed connected $M^n$ are geodesically equivalent.
 Suppose there exists a point where the metrics are strictly non-proportional.
  Then, the manifold is finitely
 covered by the product of a rational-elliptic manifold and the torus. 
\end{Cor}

\subsection{General case} 
We will sketch the   proof of  Theorem~\ref{main}. 
Consider the (1,1)-tensor  $L$ 
given by the formula (\ref{l}). 
All its eigenvalues are real. 
At every point $x\in M^n$, let  
 us denote them by  $\lambda_1(x)\le ... \le \lambda_n(x)$.
 It appears that they are globally ordered~\cite{dedicata,hyperbolic,threemanifold}: 

\begin{Th}\label{ordered123}
Let $(M^n, g)$ be a  connected Riemannian manifold.  
Suppose every  two points of the manifold can be connected by a  geodesic. 
Let Riemannian metric $\bar g$ on $M^n$ be geodesically equivalent to $g$. 

 Then, 
 for every  $i\in \{1,... ,n-1\} $, for all   $x,y\in M^n$, 
the following holds: 
\begin{enumerate}

\item $\lambda_i(x)\le \lambda_{i+1}(y)$. 

\item  If $\lambda_i(x)< \lambda_{i+1}(x)$ for some $x\in M^n$,  
then $\lambda_i(z)< \lambda_{i+1}(z)$
for almost every point $z\in M^n$.
 
\item If  $\lambda_i(x) = \lambda_{i+1}(y)$,  then there exists 
        $z\in M^n$  such that $\lambda_i(z) = \lambda_{i+1}(z)$.
\end{enumerate}
\end{Th}

Thus, if $g$ and $\bar g$ on 
closed connected $M^n$ are geodesically equivalent, then the following two cases  are possible:

\begin{itemize}                                                                         
 \item[{\bf Case 1:}] There  exists $r\in \{1,...,n-1\}$ and a constant 
  $\lambda \in \mathbb{R}$ such that, 
             for every  $x\in M^n$
$$
 \lambda_r(x)<\lambda <\lambda_{r+1}(x).
$$
   \item[{\bf Case 2:}] The following two conditions hold: 
  \begin{itemize}
 
    \item[$(i)$]For every   $r\in \{1,...,n-1\}$, the maximum 
    $\max_{x\in M^n}(\lambda_r(x))$ is equal 
    to the minimum $\min_{x\in M^n}(\lambda_{r+1}(x))$. 

    \item[$(ii)$] At least one of the eigenvalues of $L$  is not  constant.
   \end{itemize}
\end{itemize} 

In the first case,  it is possible to canonically 
construct~\cite{hyperbolic} 
the metric $\tilde g$ and the foliations $B_r$ and $B_{n-r}$ as in 
Theorem~\ref{main}. 

Let us explain where the foliations are coming from. 
At every point $x\in M^n$, let us denote by $V_r(x)$ ($V_{n-r}(x)$, respectively) 
the   sum $$ \bigoplus_{i=1}^r E_{\lambda_i} \ \ (\bigoplus_{i=r+1}^{n} E_{\lambda_i}, \ \ \textrm{respectively})$$ 
where $E_{\lambda_i}\in T_xM^n$  is the 
 eigenspace corresponding to $\lambda_i$. 

Under the assumptions of Case 1, $V_r$ and $V_{n-r}$ are smooth distributions 
of dimensions $r$ and $n-r$. By Theorem~\ref{nijenhuis}, they are integrable. Then, they generate two foliations $B_r$ and $B_{n-r}$. 

 The construction of the metric $\tilde g$ from Theorem~\ref{main}  is based on the classical Levi-Civita's  Theorem~\cite{Levi-Civita}.

In the second case, 
it is possible to    show that   the fundamental group of the manifold is 
finite.  The key  instrument for it 
is Theorem 6 of the paper~\cite{hyperbolic}  which,  
roughly speaking, tells us  that every  (closed)
 manifold with two geodesically equivalent metrics
satisfying $(i), (ii)$  has a closed 
 submanifold $U$  with two  
geodesically equivalent metrics satisfying  $(i), (ii)$ such that  
 the natural  homomorphism  $\Id_*:\pi_1(U)\to \pi_1(M^n)$
is a surjection.  
Consequently applying this theorem, 
we  come to one of the following  subcases:

\begin{itemize}
\item[{\bf SC 1:}]  The dimension  $n$ of the manifold $M^n$ 
 is $q+1$, where $q\ge 1$.  
The eigenvalues $\lambda_1=...=\lambda_q\eqdef \lambda$ are constant, the 
 eigenvalue $\lambda_{q+1}$ is not  constant and there exists 
$z\in M^{q+1}$ such that  
$\lambda_{q+1}(z)=\lambda$.

\item[{\bf SC 2:}]  The dimension $n$ of the manifold $M^n$ is $2$. 
The eigenvalues $\lambda_1$ and $\lambda_{2}$
 are  not constant and there exists a point $z\in M^{2}$ such that  
$\lambda_1(z)=\lambda_{2}(z)$.

\item[{\bf SC 3:}]  The dimension $n$ of the manifold  $M^n$ is $q+2$, 
where $q\ge 1$,  
the eigenvalues $\lambda_1$ and $\lambda_{q+2}$
 are  not constant and there exist  $z_1,z_2\in M^{n}$ such that  
$\lambda_1(z_1)=\lambda_{q+2}(z_2)$. 

\item[{\bf SC 4:}] The dimension $n$  of the manifold $M^n$  is $n=q+r+1$; 
$q>0, r>0$.  
The eigenvalues $\lambda_1=\lambda_2=...=\lambda_r$ and 
$\lambda_{r+2}=\lambda_{r+3}=...=\lambda_{n}$ are constant. 
The eigenvalue $\lambda_{r+1}$ is not constant. 
There exist 
points $z_0, z_1\in M^{n}$  such that 
$\lambda_{r+1}(z_0)=\lambda_{1}$ and 
$\lambda_{r+1}(z_1)=\lambda_{n}$.
\end{itemize}

 It is possible to  show~\cite{hyperbolic}
 that in all four subcases 
 the fundamental group in finite.   

\section*{Acknowledgments}
\noindent  
 In this paper I collected     results of  three years work; I
 am very  grateful to many   different people for 
 their interest in  this  problem.
  Especially,  I would like to thank   W. Ballmann, V. Bangert,   A. Bolsinov, 
K. Burns,   A. Fomenko, M. Gromov, U. Hamenst\"adt, M. Igarashi, 
 K. Kiyohara, B. Kruglikov,   A. Naveira,
 P. Seidel,  P. Topalov and K. Voss  for  fruitful discussions.  

  I would like to thank The University  of Tromso 
(where the results  of Section~\ref{kruglikov} were obtained)
 for the hospitality.   

 I also would like to  thank  The  
European Post-Doctoral Institute, The Max-Planck Institute for Mathematics  
(Bonn)   and   The Isaac Newton Institute for Mathematical Sciences
 for hospitality and  partial financial support.  
My research at INIMS has been  supported   by EPSRC  grant GRK99015. My research was  partially supported by 
 DFG-programm 1154 (Global Differential
Geometry)  and  Ministerium f\"ur Wissenschaft, Forschung und
Kunst  Baden-W\"urttemberg  (Elitef\"orderprogramm Postdocs 2003).

\end{document}